## On the Limit Distributions Associated with Step-kind Boundary Problems

Sherzod M. Mirakhmedov and Anatoliy N. Starscev

Institute of Mathematics. Academy of Sciences. Tashkent. Uzbekistan

**Abstract**. A random walk generated by a sum of independent identity distributed random variables with positive expectation is considered. The limiting distributions for the first- passage -time of a step-function boundary are derived.

1. **Introduction**. Let $\xi_1, \xi_2, ...$ be a sequence of i.i.d. random variables (r.v.) with

$$E\xi_i = a \text{ and } Var\xi_i = 1. \tag{1.1}$$

Set

$$S_0 = 0,\ S_n = \sum_{i=1}^{n} \xi_i,\ n = 1, 2, ...;\ \bar{S}_n = \max(S_1, ..., S_n), \tag{1.2}$$

and define in the plane $(t, \eta)$ a random polygon $\eta_n(t)$ with vertices at the points $(k/n, S_n), k = 1, ..., n$. Now we define the first-passage time of this polygon of a curve $g_n(t)$, $0 \le t \le 1$, $g_n(0) = 0$,

$$\tau_g = \min\{k : S_k \ge g_n(k/n), k = 0, 1, ..., n\} \tag{1.3}$$

If $\xi_i \ge 0$ and $g_n(t) = x_n$ for all $t \in [0,1]$, then $P\{\tau_g > n\} = P\{S_n < x_n\}$, and hence by the Central Limit Theorem we have

$$\lim_{n\to\infty} P\left\{ \frac{\tau_g - \frac{x_n}{n}}{n^{-3/2}\sqrt{x_n}} < u \right\} = \Phi(u) = \frac{1}{\sqrt{2\pi}} \int_{-\infty}^{u} e^{-z^2/2} dz, \tag{1.4}$$

uniformly in $u, |u| < \infty$. This is integral limit theorem for the number of renewals, see Feller (1966). For more general situation when just $a > 0$ the relation (1.4) follows from the fact that $P\{\tau_g > n\} = P\{\bar{S}_n < x_n\}$ and result of Wald[1] (1947) :

$$\lim_{n\to\infty} P\{\bar{S}_n < u\sqrt{n} + na\} = \Phi(u). \tag{1.5}$$

[1] Actually Wald (1947) have considered the triangular array of r.v.s $\xi_{ni}$, $n = 1, 2, ...;\ i = 1, ..., n$ and the relation (1.5) have been proved under the conditions $Var\xi_{ni} = 1$ and $\sqrt{n}E\xi_{ni} \to \infty$ as $n \to \infty$.

The first-passage time problem has been studied in the literature for the case when the boundaries are smooth enough. For the case $a=0$ the qualitative results are, in fact, the consequences of the well-known Donsker-Prokhorov's invariance principle. The first results of this kind for the case when $g_n(t)$ is a curved strip was proved by Kolmogorov and Khinchin, feature generalization of his results was done by Prokhorov, Skorokhod, Borovkov and others, for the details we refer to survey paper of Borovkov and Koroluk (1965). We refer also to S.Nagaev (1970) who have obtained the lower (w.r.t. *n*) bound for the remainder term in the limit theorem for the first passage time of a curved strip for the case $a=0$. Such kind result for the case $a>0$ have been proved by Rogozin (1966). The best bound for the remainder term in the limit theorems for the $\tau_g$ with $g_n(t)=x_n, t\in[0,1]$, is derived by Mirakhmedov (1978) for the both $a=0$ and $a>0$ cases.

The aim of the present paper is to study an asymptotic behavior of $\tau_g$ in the case when the boundary $g_n(t)$ is a piecewise constant function ( step- function) viz.,

$$g_n(t)=g_{n,i}\overset{def}{=}g_n(\tilde{t}_i),\ \tilde{t}_{i-1}<t\le\tilde{t}_i,\ \tilde{t}_l=N_l/n,\ l=1,...,k, \tag{1.6}$$

where an integer $k\ge1$, $N_l=N_l(n)$ is a sequences of integers numbers such that

$$0=N_0(n)\le N_1(n)\le...\le N_k(n)=n \text{ and } \lim_{n\to\infty}N_l(n)=\infty.$$

**Remark 1.1**. If the possible values of the r.v. $\xi_k$ are $-1$ and $+1$ only, then $\tau_g$ can be interpreted as the moment of bankruptcy of a player who adds to his pocket amount of $g_n(\tilde{t}_l)-g_n(\tilde{t}_{l-1})$ at the moments $N_l$ in a game with an infinitely rich partner.

We also note that the functional $\tau_g$ is not continuous due to discontinuity of the boundary $g_n(t)$, and hence it is impossible to use Donsker-Prokhorov's invariance principle.

2. **Results**. In what follows we deal with the case $a>0$, also to keep notations simple we shall assume that $g_n(t)$ is non-decreasing function. Due to (1.2), (1.3) and (1.6) one can observe that

$$P\left\{\tau_g>n\right\}=P\left\{\bar{S}_l<g_n(\tilde{t}_l),\ l=1,...,k-1,\ \bar{S}_n<g_n(1)\right\}. \tag{2.1}$$

Set $\bar{S}^*_{N_l}=\left(\bar{S}_{N_l}-N_l(n)a\right)/\sqrt{N_l(n)}$.

**Theorem 1**. Let $\xi_1,\xi_2,...$ be i.i.d. r.v.s with $E\xi_1=a>0$ and $Var\xi_i=1$. Then

$$\lim_{n\to\infty}\sup_{u_1,...,u_k}\left[P\left\{\bar{S}^*_{N_1}<u_1,...,\bar{S}^*_{N_k}<u_k\right\}-\Phi_{\Lambda_k}(u_1,...,u_k)\right]=0,$$

where $\Phi_{\Lambda_k}(.)$ is a *k*-dimensional normal distribution function with zero vector of expectations and covariance matrix $\Lambda_k=\left(\sqrt{\lambda_{ij}}\right)$, $\lambda_{ij}=N_i(n)/N_j(n), 1\le i\le j\le k$.

**Remark**. The quantities $\lambda_{ij}$ are defined for $1\le i\le j\le k$, for $i>j$ we have to put $\lambda_{ij}$ as $\lambda_{ij}$ with $i<j$.

In particular, if $k=1$ and $N_1(n)=n$ then Theorem 1 implies relation (1.5).

Now let us turn to the description of the asymptotic distribution of $\tau_g$. Set, see notation (1.6),

$$V_{n,i}=\frac{N_i(n)-g_{n,i}a^{-1}}{a^{-3/2}\sqrt{g_{n,i}}}\quad,\quad \Delta_{n,i}=\frac{g_{n,i+1}-g_{n,i}}{\sqrt{g_{n,i+1}/a}},\quad \lim_{n\to\infty}V_{n,i}=V_i,\quad \lim_{n\to\infty}\Delta_{n,i}=\Delta_i. \qquad (2.2)$$

We shall assume that $V_{i-1}<V_i, i=1,...,k$; otherwise at the limit no uniqueness of the normalized boundary will appear (due to the influx of its step-parts) and the situation becomes to be simpler.

Set $\tau_g^{(i)}=\left(\tau_g-g_{n,i}a^{-1}\right)/a^{-3/2}\sqrt{g_{n,i}}$.

**Theorem 2**. Let $\xi_1,\xi_2,...$ be i.i.d. r.v.s satisfying condition (1.1) with $a>0$ and

$$N_i(n)/N_j(n)=\lambda_{ij}(1+o(1)),\ N_i(n)=g_{n,i}/a\left(1+o(1)\right),\ i=1,...,k_0<k \text{ and } V_{k_0+1}=\infty, \qquad (2.3)$$

Then for each $i=1,...,k_0+1$ and arbitrary fixed $u$ one has

$$\lim_{n\to\infty}P\left\{\tau_g^{(i)}<u\right\}=G_i(u),$$

where

$$G_i(u)=\begin{cases}1-\Phi_{\Lambda_j}\left(-V_1,...,-V_{j-1},-u-\Sigma_j^{i-1}\right)\ if\ \ \lambda_{j-1,i}^{1/2}V_{j-1}-\Sigma_{j-1}^{i-1}<u\le\lambda_{j,i}^{1/2}V_j-\Sigma_j^{i-1}, 1\le j\le i,\\ 1-\Phi_{\Lambda_j}\left(-V_1,...,-V_{j-1},-u+\Sigma_j^{i-1}\right)\ if\ \ \lambda_{j-1,i}^{1/2}V_{j-1}+\Sigma_i^{i-2}<u\le\lambda_{j,i}^{1/2}V_j+\Sigma_i^{j-1}, i<j\le k_0+1\end{cases},$$

here we put $V_0=-\infty$, $\lambda_{0,i}=1$, $\Sigma_i^j=0$ for $j\le i$ and $\Sigma_i^j=\sum_{s=1}^{i}\Delta_s\lambda_{s+1,i}^{1/2}$ for $j>i$.

Note that if

$$|V_i|<\infty, \Delta_i<\infty,\ i=1,...,k_0,\ V_{k_0+1}=\infty, \qquad (2.4)$$

then all $G_i$ are proper distributions, and hence each of them well enough characterizes the r.v. $\tau_g$. In this case the limit distributions $G_i$ are simplifies, since $\lambda_{rs}=1$ for $1\le r\le s$ due to finiteness of $\Delta_i$ s.

**Corollary 1**. Under the condition (2.4) we have

$$G_i(u)=\begin{cases}\Phi(u)\ \text{ if }\ \ u\le V_i,\\ \Phi(V_i)\ \text{ if } V_i+\sum_{j=1}^{i-1}\Delta_j<u\le V_i+\sum_{j=1}^{i}\Delta_j, i=1,...,k_0,\\ \Phi\left(u-\sum_{j=1}^{i}\Delta_j\right)\ \text{ if } V_i+\sum_{j=1}^{i}\Delta_j<u\le V_{i+1}+\sum_{j=1}^{i}\Delta_j, i=1,...,k_0-1,\\ \Phi\left(u-\sum_{j=1}^{k_0}\Delta_j\right)\ \text{ if } u>V_{k_0}+\sum_{j=1}^{k_0}\Delta_j.\end{cases}$$

In particular, if $\Delta_i=0$ for $i=1,...,k_0$ then $G_1(u)=\Phi(u)$; if $V_1=\infty$ then $G_1(u)=\Phi(u)$ also.

Note that under the conditions (2.3)

$$\lim_{n\to\infty} P\left\{V_i + \sum_{j=1}^{i-1}\Delta_j < \tau_g^{(1)} \le V_i + \sum_{j=1}^{i}\Delta_j\right\} = 0,\ i = 1,...,k_0.$$

If $\Delta_i = \infty$ for $i = 1,...,k_0$ then $G_i$ s becomes non-proper distributions. For example, if $\Delta_1 = \infty$ then it follows from Theorem 2 that $G_1(u) = \Phi(\min(u,V_1))$.

**Corollary 2**. If $\Delta_i = \infty$, $\lim_{n\to\infty} g_{n,i} / g_{n,i+1} = \alpha_i$, $i = 1,...,k_0$, and $V_{k_0} = \infty$, then

$$G_i(u) = 1 - \Phi_{\Lambda_j^*}\left(-V_1,...,-V_{j-1}, -\min(u,V_i)\right) \text{ for } i = 1,...,k_0$$

and

$$G_{k_0+1}(u) = 1 - \Phi_{\Lambda_{k_0+1}^*}\left(-V_1,...,-V_{j-1}, -u\right),$$

where $\Lambda_i^*$ is a symmetric matrix with elements $\lambda_{rs}^{1/2} = \prod_{j=r}^{s} \alpha_j^{1/2}$, $r < s \le i$, $\lambda_{rr} = 1$.

It is easy to see that all distributions $G_i$ in Corollary 2 are non-proper, but they are complemented of each other such that a sum of their total variations is equal to 1. Indeed,

$$\sum_{i=1}^{k_0+1}\left[G_i(+\infty) - G_i(-\infty)\right] = 1 - \Phi(-V_1)$$

$$+\sum_{i=2}^{k_0}\left[\Phi_{\Lambda_{i-1}^*}(-V_1,...,-V_{i-1}) - \Phi_{\Lambda_i^*}(-V_1,...,-V_i)\right] + \Phi_{\Lambda_{k_0}^*}(-V_1,...,-V_{i-1}) = 1.$$

**Corollary 3**. If $\Delta_i = \infty$, $\lim_{n\to\infty} g_{n,i} / g_{n,i+1} = \alpha_i$, $i = 1,...,k_0$, $V_{k_0} = \infty$ and $\alpha_i = 1$, $i = 1,...,k_0$ then

$$G_i(u) = \begin{cases} \Phi(V_{i-1}) & \text{if} \quad u \le V_{i-1}, \\ \Phi(u_i) & \text{if} \quad V_{i-1} < u \le V_i, \\ \Phi(V_i) & \text{if} \quad u > V_i. \end{cases}$$

$$G_{k_0+1}(u) = \begin{cases} \Phi(V_{k_0}) & \text{if} \quad u \le V_{k_0}, \\ \Phi(u) & \text{if} \quad u > V_{k_0}. \end{cases}$$

We are restricted ourselves by these corollaries, while it is clear that under another combination of conditions on $V_i$ and $\Delta_i$ one can derive corresponding consequences of Theorem 2. Instead we shall consider the rate of convergence in Theorem 1, such result can be applied to study the rate of convergences in Theorem 2 and its corollaries.

**Theorem 3**. Let $\xi_1, \xi_2, ...$ be i.i.d. r.v.s such that $E\xi_i = a > 0$, $Var\xi_i = 1$ and $\beta_3 = E\left|\xi_i - a\right|^3$. Then for arbitrary $1 \le N_1(n) \le ... \le N_k(n)$ the following is hold

$$\sup_{u_1,...,u_k}\left|P\left\{\bar{S}^*_{N_1}<u_1,...,\bar{S}^*_{N_k}<u_k\right\}-\Phi_{\Lambda_k}(u_1,...,u_k)\right|\le\tilde{C}_1\sum_{i=1}^{k}\frac{1}{\sqrt{N_i(n)}}+C_0\beta_3\sum_{i=1}^{k}\frac{2^{k-1}}{\sqrt{N_i(n)-N_{i-1}(n)}},$$

where $C_0 \le 0.82$, $\tilde{C}_1 = C_1\left(\max\left(E|\xi_1|,\beta_3,\beta_3/a\right)\right)^2$, the universal constant $C_1$ is from Theorem of Rogozin (1966), and $\Lambda_k$ is a symmetric matrix with elements $\lambda_{ij}=\sqrt{N_i(n)/N_j(n)}$, for $1\le i\le j\le k$.

3. **Proofs**. We shall use the notations of the previous sections.

**Proof of Theorem 1**. First we show that for arbitrary $u_1,...,u_k$

$$\delta_{n,k}=P\left\{S_{N_1}<u_1,...,S_{N_k}<u_k\right\}-P\left\{\bar{S}_{N_1}<u_1,...,\bar{S}_{N_k}<u_k\right\}\to 0 \text{ as } n\to\infty. \qquad (3.1)$$

We shall proof (3.1) by induction w.r.t. $k$. For $k=1$ relation (3.1) follows from (1.5) and CLT for i.i.d. r.v.s . Assume that (3.1) is hold true for $k=j$ ( i.e. $\delta_{n,j}\to 0$) and we wish to prove it for $k=j+1$. We have $\delta_{n,j+1}=\delta_{n,j}+\delta'_{n,j}$, where

$\delta'_{n,j}=P\left\{\bar{S}_{N_1}<u_1,...,\bar{S}_{N_j}<u_j,\bar{S}_{N_{j+1}}\ge u_{j+1}\right\}-P\left\{S_{N_1}<u_1,...,S_{N_j}<u_j,S_{N_{j+1}}\ge u_{j+1}\right\}$. If $\delta'_{n,j}\le 0$ then $\delta_{n,j+1}\le\delta_{n,j}\to 0$. Let $\delta'_{n,j}>0$, then

$$\begin{aligned}\delta'_{n,j}&\le P\left\{S_{N_1}<u_1,...,S_{N_j}<u_j,\bar{S}_{N_{j+1}}\ge u_{j+1}\right\}-P\left\{S_{N_1}<u_1,...,S_{N_j}<u_j,S_{N_{j+1}}\ge u_{j+1}\right\}\\ &\le P\left\{S_{N_1}<u_1,...,S_{N_j}<u_j,\bar{S}_{N_{j+1}}\ge u_{j+1},S_{N_{j+1}}<u_{j+1}\right\}\le P\left\{\bar{S}_{N_{j+1}}\ge u_{j+1},S_{N_{j+1}}<u_{j+1}\right\}\\ &=P\left\{S_{N_{j+1}}<u_{j+1}\right\}-P\left\{\bar{S}_{N_{j+1}}<u_{j+1}\right\}\to 0.\end{aligned}$$

Thus $\delta_{n,j+1}\to 0$ and ( 3.1) follows by induction. Note that (3.1) implies

$$P\left\{S^*_{N_1}<u_1,...,S^*_{N_k}<u_k\right\}-P\left\{\bar{S}^*_{N_1}<u_1,...,\bar{S}^*_{N_k}<u_k\right\}\to 0, \qquad (3.2)$$

as $n\to\infty$ uniformly in $u_1,...,u_k$, where $\bar{S}^*_{N_i}=\left(\bar{S}_{N_l}-N_l(n)a\right)/\sqrt{N_l(n)}$ and $S^*_{N_i}=\left(S_{N_l}-N_l(n)a\right)/\sqrt{N_l(n)}$. On the other hand by CLT for i.i.d. random vectors in $R^k$, see for instance Feller (1966, Chapter VIII, Theorem 2) $P\left\{S^*_{N_1}<u_1,...,S^*_{N_k}<u_k\right\}-\Phi_{\Lambda_k}(u_1,...,u_k)\to 0$ as $n\to\infty$ uniformly in $u_1,...,u_k$. Theorem 1 follows from this fact and (3.2).

**Proof of Theorem 2**. Consider a sequence $\{k_n\}$ such that $N_{j-1}(n)<k_n\le N_j(n)$. Then due to (2.1) we have

$$P\left\{\tau_g>k_n\right\}=P\left\{\bar{S}_l<g_{n,l},\, l=1,...,j-1,\, \bar{S}_{k_n}<g_{n,j}\right\}. \qquad (3.3)$$

Let $j\le i$. Due to (3.3) we have, see notations (2.2) and (1.6),

$$P\left\{\tau_g^{(i)} > \left(k_n - g_{n,i}a^{-1}\right)/ a^{-3/2}\sqrt{g_{n,i}}\right\}$$

$$= P\left\{\bar{S}^*_{N_1} < -V_{n,1}\sqrt{\frac{g_{n,1}}{N_1 a}}, ..., \bar{S}^*_{N_{j-1}} < -V_{n,j-1}\sqrt{\frac{g_{n,j-1}}{N_{j-1}a}},\ \bar{S}^*_{N_{k_n}} < -\left(\frac{k_n - g_{n,i}a^{-1}}{a^{-3/2}\sqrt{g_{n,i}}} + \frac{g_{n,i}-g_{n,j}}{\sqrt{g_{n,i}/a}}\right)\sqrt{\frac{g_{n,i}}{k_n a}}\right\}. \quad (3.4)$$

Next, due to condition (2.3) we have $\lim_{n\to\infty} g_{n,i} / g_{n.j} = \lambda_{ij}$, and hence

$$\lim_{n\to\infty}\frac{g_{n,i}-g_{n,j}}{\sqrt{g_{n,i}/a}} = \lim_{n\to\infty}\sum_{s=j}^{i-1}\Delta_{n,s}\sqrt{\frac{g_{n,s+1}}{g_{n,i}}} = \sum_{s=j}^{i-1}\Delta_s\sqrt{\lambda_{s+1,i}} \stackrel{def}{=} \Sigma_j^{i-1}. \quad (3.5)$$

Choosing now $k_n$ such that $\lim_{n\to\infty}\left(k_n - g_{n,i}a^{-1}\right)/ a^{-3/2}\sqrt{g_{n,i}} = u$ form (3.4) , (3.5) and Theorem 1 we obtain

$$\lim_{n\to\infty} P\left\{\tau_g^{(i)} > u\right\} = \Phi_{\Lambda_j}\left(-V_1, ..., -V_{j-1}, -u - \Sigma_j^{i-1}\right). \quad (3.6)$$

Because $k_n \in (N_{j-1}, N_j]$ the range of variable *u* in (3.6) is

$$\lim_{n\to\infty}\left(V_{n,j-1}\sqrt{\frac{g_{n,j-1}}{g_{n,i}}} - \frac{g_{n,i}-g_{n,j-1}}{\sqrt{g_{n,i}/a}}\right) < u = \lim_{n\to\infty}\frac{k_n - g_{n,i}a^{-1}}{a^{-3/2}\sqrt{g_{n,i}}} \le \lim_{n\to\infty}\left(V_{n,j}\sqrt{\frac{g_{n,j}}{g_{n,i}}} - \frac{g_{n,i}-g_{n,j}}{\sqrt{g_{n,i}/a}}\right).$$

This implies $V_{j-1}\sqrt{\lambda_{j-1,i}} - \Sigma_{j-1}^{i-1} < u \le V_j\sqrt{\lambda_{j,i}} - \Sigma_j^{i-1}$, since (2.3) and (3.5).

The case $j > i$ can be studied analogically.

**Proof of Theorem 3**. First we note that for arbitrary $u_1, ..., u_k$, see (3.1),

$$\delta_{n,k} \le C_1\sum_{i=1}^{k} N_i^{-1/2}. \quad (3.7)$$

Indeed, for $k = 1$ (3.7) is the result of Rogozin (1966). For $k > 1$, the inequality (3.7) can be proved by induction, reasoning analogically to proof of (3.1).

Now it is remain to prove that

$$J_k = \sup_{u_1,...,u_k}\left|P\left\{S^*_{N_1} < u_1, ..., S^*_{N_k} < u_k\right\} - \Phi_{\Lambda_k}(u_1, ..., u_k)\right| \le C_0\beta_3\sum_{j=1}^{k}\frac{2^{k-j}}{\sqrt{N_j - N_{j-1}}}. \quad (3.8)$$

Note that this inequality (3.8) can be derived from Berry-Esseen bound result in k-dimensional CLT, see for instance Corollary 15.3 of Bhattacharya and Rao (1976), but then we need to calculate the determinant and minors of the matrix $\Lambda_k$, and hence this approach need complex algebra. Instead we shall use again induction in *k* method.

For k =1 inequality (3.8) is well-known Theorem of Berry-Esseen, see Feller (1966, Ch. XVI, Sec,5). For k=2 we write

$$J_2 \le \sup_{u_1,u_2}\left|\int_{-\infty}^{u_1}\Phi\left(\frac{u_2 - \lambda_{12}z}{\sqrt{1-\lambda_{12}^2}}\right)d\Phi(z) - \Phi_{\Lambda_2}(u_1, u_2)\right|$$

$$+\sup_{u_1,u_2}\left|\int_{-\infty}^{u_1}\left[P\left\{S^*_{N_2-N_1}<\frac{u_2-\lambda_{12}z}{\sqrt{1-\lambda_{12}^2}}\right\}-\Phi\left(\frac{u_2-\lambda_{12}z}{\sqrt{1-\lambda_{12}^2}}\right)\right]dP\left\{S^*_{N_1}<z\right\}\right|$$

$$+\sup_{u_1,u_2}\left|\int_{-\infty}^{u_1}\Phi\left(\frac{u_2-\lambda_{12}z}{\sqrt{1-\lambda_{12}^2}}\right)d\left(P\left\{S^*_{N_1}<z\right\}-\Phi(z)\right)\right|=J_{21}+J_{22}+J_{23}.$$

By simple calculation we can see that $J_{21}=0$. Using Berry-Esseen's theorem we obtain $J_{21}\le C_0\beta_3\left(N_2-N_1\right)^{-1/2}$ and integrating by part

$$J_{23}\le\sup_{u_1,u_2}\left|P\left\{S^*_{N_1}<z\right\}-\Phi(z)\right|\Phi\left(\frac{u_2-\lambda_{12}z}{\sqrt{1-\lambda_{12}^2}}\right)$$

$$+\sup_{u_1,u_2}\int_{-\infty}^{u_1}\left|P\left\{S^*_{N_1}<z\right\}-\Phi(z)\right|d\Phi\left(\frac{u_2-\lambda_{12}z}{\sqrt{1-\lambda_{12}^2}}\right)\le 2C_0\beta_3N_1^{-1/2}.$$

So

$$J_2\le\frac{\le 2C_0\beta_3}{\sqrt{N_1}}\left(1+\frac{\lambda_{12}}{2\sqrt{1-\lambda_{12}^2}}\right).$$

Analogically for any $j>1$ we can write

$$P\left\{S^*_{N_1}<u_1,\dots,S^*_{N_j}<u_j\right\}=\int_{-\infty}^{u_{j-1}}P\left\{S^*_{N_j-N_{j-1}}<\frac{u_j-\lambda_{j-1,j}z_{j-1}}{\sqrt{1-\lambda_{j-1,j}^2}}\right\}d_{z_{j-1}}P\left\{S^*_{N_1}<u_1,\dots,S^*_{N_{j-1}}<z_{j-1}\right\}.$$

Therefore, reasoning alike to the case $k=2$ by induction w.r.t. k we can complete the proof of (3.8). Theorem 3 is proved.